\documentclass{amsart}
\usepackage{amsthm,epsfig}
\usepackage{amsmath}
\usepackage{amssymb, amstext, epsfig}
\theoremstyle{plain}
\newtheorem{lemma}{Lemma}[section]

\newtheorem{theorem}[lemma]{Theorem}

\theoremstyle{remark}
\newtheorem{remark}[lemma]{Remark}

\newcommand{\Z}{{\mathbb Z}}

\newcommand{\R}{{\mathbb R}}
\newcommand{\C}{{\mathbb C}}

\newcommand{\eps}{{\varepsilon}}

\title{
Linear estimate for the number of zeros of Abelian integrals}

\author{Sergey Malev} \address{Faculty of Mathematics and Computer
Science, The Weizmann Institute of Science POB 26, Rehovot 76100,
ISRAEL}\email{sergey.malev@weizmann.ac.il}

\author{Dmitry Novikov} \address{Faculty of Mathematics and Computer
Science, The Weizmann Institute of Science POB 26, Rehovot 76100,
ISRAEL}\email{dmitry.novikov@weizmann.ac.il}

\begin{document}

\begin{abstract}
We prove a linear in $\deg\omega$ upper bound on the number of real zeros of the Abelian integral $I(t)=\int_{\delta(t)}\omega$, where $\delta(t)\subset\R^2$ is the real oval $x^2y(1-x-y)=t$ and $\omega$ is a one-form with polynomial coefficients.
\end{abstract}
\maketitle
\section{Introduction}

For the polynomial $H(x,y)=x^2y(1-x-y)$ consider the  continuous family
$\{\delta(t), t\in\left(0,\frac{1}{64}\right)\}$ of compact connected components of the level curves $\{H=t\}\subset\R^2$. Let $\omega=p(x,y)dx+q(x,y)dy\in\Lambda^1(\R^2)$ be a differential one-form with polynomial coefficients of degree $n$. Define the complete Abelian integral:
\begin{equation}\label{def_abelint}
I(t)=\oint_{\delta(t)}\omega,\ \ \ t\in\left(0,\frac{1}{64}\right).
\end{equation}

We provide an explicit answer to the Infinitesimal Hilbert 16th problem for this particular Abelian integral.
\begin{theorem}\label{main_theorem}
The number of isolated zeros of $I(t)$ on $(0,\frac{1}{64})$
 does not exceed $\frac{7}{4}n+9$, where $n=\deg \omega$.
\end{theorem}

It is well known that zeros of Abelian integrals correspond to
limit cycles appearing in non-conservative perturbations of Hamiltonian, or, more general, integrable systems. Abelian integral \eqref{def_abelint} is related to  perturbations of the integrable quadratic vector field which can be written in the Pfaffian form as follows:
\begin{equation}\label{pfaffian2}
\frac{1}{x}dH-\eps\omega=0, \qquad H(x,y)=x^2y(1-x-y).
\end{equation}

Therefore Theorem~\ref{main_theorem} implies in a standard way the following
claim:
\begin{theorem}\label{main_theorem:cycles}
The number of limit cycles appearing in non-conservative perturbation \eqref{pfaffian2} and converging, as $\epsilon\to0$, to a smooth cycle $\delta(t)$,  does not exceed $\frac{1}{4}(7n+43)$.
\end{theorem}
Indeed, these limit cycles correspond to the isolated zeros of $\int_{\delta(t)}x\omega$, so this upper bound follows from Theorem~\ref{main_theorem} by replacing $n$ by $n+1$.

Our result should be considered in the general context of the Infinitesimal Hilbert 16th problem. So far, the only known general explicit result about the number of zeros of Abelian integrals is the recent result \cite{BNY} providing double-exponential in $\max(\deg H,\deg \omega)$ upper bound for the number of zeros of Abelian integral. The result of Petrov-Khovanskii, see \cite{Z} for an exposition of the result,  provides an  upper bound which is a linear function of  $n=\deg\omega$, but provides no information about the coefficients of this function. It seems reasonable to expect that these two results can be combined together to provide an upper bound which would be  linear in $n$ and double-exponential in $\deg H$. However, even this upper bound will by far exceed any known examples.

From the other side, for the cases of polynomials $H$ of low degree the situation seems to be better understood. More exact, for a generic polynomial $H$ of third degree Horozov and Iliev \cite{HI} were able to provide an explicit upper bound linear in $n$ which seems to be close to the best possible one. The key point was  the ellipticity of the level curves of $H$. This fact allowed to reduce the initial question to the question about number of zeros of polynomial combinations of solutions of some Riccati equation. This last question can be easily dealt with by essentially fewnomials technique (i.e. Rolle lemma). Later, the same approach was applied to the case of elliptic polynomial of fourth degree, see \cite{GJ}.

This result motivated consideration of integrable quadratic vector fields with center whose trajectories are elliptic curves. Gautier in \cite{Gau} lists all these fields. Based on this list and results of \cite{I}, Gautier, Gavrilov and Iliev in \cite{GGI} proposed a program of studying cyclicity of open nests of cycles defined by such foliations, and, in particular,  conjectured that one can provide an effective upper bound for the number of zeros of corresponding Abelian integrals, similar to  \cite{HI}. Our case is the case (rlv3) in notations of \cite{GGI}.

\section{Decomposition in Petrov modules}
Define the basic Abelian integrals as
\begin{equation}\label{basic_abel}
I_{i,j}(t)=\frac{1}{i+1}\int_{\delta(t)}x^{i+1}y^jdy=\iint_{\Delta(t)}x^iy^jdx\wedge dy,
\end{equation}
where $\Delta(t)$ is the area bounded by the cycle $\delta(t)$. The integrals are well-defined for all  (i.e. not only positive) $i,j\in\Z$ due to the following fact:

\begin{remark}\label{lem_pos}
The curve $\delta(t)$ lies in $\{x,y>0\}$.
\end{remark}

Our immediate goal is to construct explicit representation of the Abelian integrals defined in (\ref{def_abelint}) as combinations with polynomial in $t$ coefficients of just three Abelain integrals $J_1(t)=I_{0,0}(t)$, $J_2(t)=I_{2,0}(t)$ and
$J_3(t)=I_{3,0}(t)$. In other words, we want to prove that Abelian integrals can be generated, as a $\C[t]$-module, by these 3 basic Abelian integrals.

Let $H(x,y)=\sum_{i,j\in\Z^2}h_{ij}x^iy^j\in\R[x,y]$ be a general Laurent polynomial in two variables and assume that the family $\delta(t)\subset\{H=t\}$ of cycles lies in $\{x,y>0\}$.
\begin{lemma}\label{lemma}
Abelian integrals $I_{k,l}(t)$ defined by \eqref{basic_abel} satisfy the following relations:
\begin{align}
t(k+1)\cdot I_{k,l}(t)&=\sum\limits_{i,j}h_{i,j}(k+i+1)\cdot I_{k+i,l+j}(t),\notag
\\
\label{abel_sys_gen}
t(l+1)\cdot I_{k,l}(t)&=\sum\limits_{i,j}h_{i,j}(l+j+1)\cdot I_{k+i,l+j}(t),\\
(k+1)\cdot I_{k,l}(t)&=\sum\limits_{i,j}ih_{i,j}\cdot\frac{d}{d t}I_{k+i,l+j}(t),\notag\\
(l+1)\cdot I_{k,l}(t)&=\sum\limits_{i,j}jh_{i,j}\cdot\frac{d}{d t}I_{k+i,l+j}(t),\notag
\end{align}
where $i,j\in\Z$.
\end{lemma}
\begin{proof}
Since $H(x,y)=t$ on the cycle $\delta(t)\subset\R^2$ we have
$$t(k+1)\cdot I_{k,l}(t)=\int\limits_{\delta(t)}H(x,y)x^{k+1}y^ldy=\sum_{i,j}h_{i,j}(k+i+1)I_{k+i,l+j}(t),$$
which is the first equality of (\ref{abel_sys_gen}).
The second identity is proved similarly.

By Gelfand-Leray formula we have
\begin{equation}\label{Gelfand-Leray}
\frac{d}{dt}(\iint\limits_{\gamma(t)}dH\wedge x^{k+1}y^ldy) =\oint\limits_{\delta(t)}x^{k+1}y^ldy=(k+1)I_{k,l}(t).
\end{equation}
Replacing $H$ by $\sum\limits_{(i,j)\in\Z^2}h_{ij}x^iy^j$ we get the
third identity. The fourth equality is proved similarly.
\end{proof}

For our particular choice $H(x,y)=x^2y(1-x-y)$ we can rewrite the relations of Lemma~\ref{lemma} in more convenient form.
We have
\begin{equation}\label{representaion_0}
\left(\begin{array}{cccc}
  t(k+1) & -k-3 & k+4 & k+3 \\
  t(l+1) & -l-2 & l+2 & l+3
\end{array}\right)
\cdot
\left(\begin{array}{c}
  I_{k,l} \\
  I_{k+2,l+1} \\
  I_{k+3,l+1} \\
  I_{k+2,l+2}
\end{array}\right)
=
\begin{pmatrix}
  0\\
  0
\end{pmatrix}  .
\end{equation}

Multiplying by
$\begin{pmatrix}
  l+3& -k-3\\
  -l-2& k+4
\end{pmatrix}$
from the left we get an equivalent for $k+l\neq -6$ system of equations:

\begin{equation}\label{representaion_1}
\left(\begin{array}{cccc}
  -2t(l-k) & -k-3 & k+l+6 & 0 \\
  t(3l-k+2) & -l-2 & 0 & k+l+6
\end{array}\right)
\cdot
\left(\begin{array}{c}
  I_{k,l} \\
  I_{k+2,l+1} \\
  I_{k+3,l+1} \\
  I_{k+2,l+2}
\end{array}\right)
=
\begin{pmatrix}
  0\\
  0
\end{pmatrix}.
\end{equation}

In other words, $I_{k+3,l+1}(t)$ and $I_{k+2,l+2}(t)$ can be represented as $\mu^{k,l}_1(t)I_{k,l}(t)+\mu^{k,l}_2(t)I_{k+2,l+1}(t)$, where $\mu^{k,l}_1$ is a
polynomial of degree at most $1$ and $\mu^{k,l}_2$ is a constant. We have such representation for any pair $(k,l)$ such that $k+l\neq -6.$

Another corollary of (\ref{representaion_0}):
\begin{equation}\label{representaion_2}
\left(\begin{array}{cccc}
 -2l+k-1 & 3l-k+2 & 2l-2k
\end{array}\right)
\cdot
\left(\begin{array}{c}
  I_{k+2,l+1} \\
  I_{k+3,l+1} \\
  I_{k+2,l+2}
\end{array}\right)
=  0.
\end{equation}
The equation (\ref{representaion_2}) gives us linear dependence between $I_{k+2,l+1},I_{k+3,l+1}$ and $I_{k+2,l+2}$ in all cases except $k=l=-1$.
But in this case the equation \eqref{representaion_1} becomes:

\begin{equation}\label{k=l=-1}
\left(\begin{array}{cccc}
  0 & -2 & 4 & 0 \\
  0 & -1 & 0 & 4
\end{array}\right)
\cdot
\left(\begin{array}{c}
  I_{-1,-1} \\
  I_{1,0} \\
  I_{2,0} \\
  I_{1,1}
\end{array}\right)
=
\begin{pmatrix}
  0\\
  0
\end{pmatrix}.
\end{equation}

Recall that $J_1(t)=I_{0,0}(t)$, $J_2(t)=I_{2,0}(t)$ and $J_3(t)=I_{3,0}(t)$.
\begin{lemma}\label{lemma32}
For any polynomial $1$-differential form $\omega$ the Abelian
integral $I(t)=\int_{\delta(t)}\omega$ can be represented as
$p_1(t)J_1(t)+p_2(t)J_2(t)+p_3(t)J_3(t)$ for some polynomials
$p_i(t)$ of degree less than or equal to $\frac{n}{4}$, where
$n=\deg \omega$.
\end{lemma}

\begin{remark} The proof of this result below essentially provides the coefficients $p_i(t)$, i.e. provides an effective decomposition in the Petrov module corresponding to $H$. This result does not formally follow from the result of Gavrilov~\cite{G} since $x^2y(1-x-y)$ is not a semi-weighted homogeneous polynomial. However, in this simple situation it can be obtained by a straightforward computation.
\end{remark}

\begin{proof}
First of all we will give such a representation for all $I_{k,l}(t)$
if $k+l\leq 3$, $k,l\geq 0$. By \eqref{k=l=-1} we have $I_{1,0}=2J_2$ and
$I_{1,1}=\frac{1}{2}J_2$. Using \eqref{representaion_2} one can
calculate the required representation for $I_{0,1}$ and $I_{2,1}$.
This implies similar representation for $I_{0,2}$ and $I_{1,2}$ and
then for  $I_{0,3}$. Note that for these integrals the coefficients
$p_i$ are polynomials of degree $0$, i.e. scalar.

For $n=i+j\ge 4$ the proof goes by induction on $n$. Using \eqref{representaion_1}, we see that the required representation of $I_{k,n-k}$  for  $2\leq k\leq n-1$, together with bounds on the degrees of $p_i$, follows from the same for $I_{k,l}$ with smaller $k+l$.
Using \eqref{representaion_2} one can obtain that
$I_{1,n-1}$ is a linear combination of $I_{1,n-2}$ and $I_{2,n-2}$ (for $k=-1$ and $l=n-3$)
because $2l-2k= 2n-4>0$. Thus we obtain a required representation for $I_{1,n-1}$ and similarly for $I_{0,n}$.
Now we use \eqref{representaion_2} for $k=n-3$ and $l=-1$ to obtain representation of $I_{n,0}$ as a linear combination of
$I_{n-1,0}$ and $I_{n-1,1}$. One can easily check that the degrees of $p_i$ are bounded by $n/4$ in all these  cases.
\end{proof}

\section{Construction of the Picard-Fuchs system.}

It is well known that the effective decomposition in Petrov modules allows to explicitly construct the Picard-Fuchs system for the generators of the Petrov module, see e.g. \cite{Y}. Here we follow this classical path.

\begin{lemma}\label{lemma_PF}
The column
$J=\begin{pmatrix}
  J_{1}(t) \\
  J_{2}(t) \\
  J_{3}(t)
\end{pmatrix}$
satisfies the system
\begin{equation}\label{Picard_Fuchs}
J=\frac{d}{dt}((A+tB)J),
\end{equation}
where $A=\begin{pmatrix}
  0 & -\frac{1}{12} & \frac{1}{12} \\
  0 & -\frac{1}{56} & \frac{1}{56} \\
  0 & -\frac{5}{504} & \frac{5}{504}
  \end{pmatrix}  $
and $B=\begin{pmatrix}
  \frac{2}{3} & 0 & 0 \\
  0 & \frac{4}{7} & 0 \\
  0 & \frac{2}{21} & \frac{4}{9}
\end{pmatrix}. $
\end{lemma}

\begin{proof}of the lemma \ref{lemma_PF}:

By Lemma \ref{lemma} for any $k$ and $l$ (and in particular for $l=0$ and $k=0,1,3$) we have
$(l+1)\cdot I_{k,l}(t)=\sum_{i,j}jh_{i,j}\cdot\frac{d}{d t}I_{k+i,l+j}(t).$

It implies
\begin{equation}\label{eq}
I_{k,0}=\frac{d}{dt}(I_{k+2,1}-I_{k+3,1}-2I_{k+2,2}).
\end{equation}
Using Lemma \ref{lemma32} we represent
$I_{k,1}$ for $2\leq k\leq 6$ and $I_{k,2}$ for $k=2, 4$ and $5$ in terms of $J_i$.
After calculation we obtain the result:
\begin{align}
I_{2,1}&=\frac{1}{2}J_2-\frac{1}{2}J_3,\notag\\
I_{3,1}&=\frac{1}{4}J_2-\frac{1}{4}J_3,\notag\\
I_{4,1}&=\frac{1}{7}J_2-\frac{1}{7}J_3-\frac{4}{7}tJ_2,\notag\\
\label{diff_sys_eq}
I_{5,1}&=\frac{5}{56}J_2-\frac{5}{56}J_3-\frac{6}{7}tJ_2,\\
I_{6,1}&=\frac{5}{84}J_2-\frac{5}{84}J_3-\frac{4}{7}tJ_2-\frac{2}{3}tJ_3,\notag\\
I_{2,2}&=\frac{1}{6}J_2-\frac{1}{6}J_3-\frac{1}{3}tJ_1,\notag\\
I_{4,2}&=\frac{1}{28}J_2-\frac{1}{28}J_3-\frac{1}{7}tJ_2,\notag\\
I_{5,2}&=\frac{5}{252}J_2-\frac{5}{252}J_3-\frac{4}{21}tJ_2+\frac{1}{9}tJ_3,\notag
\end{align}

The  formulas \eqref{eq} and system \eqref{diff_sys_eq} together immediately imply the system \eqref{Picard_Fuchs}.
\end{proof}
The system \eqref{Picard_Fuchs} can be rewritten as
\begin{equation}
(A+tB)\cdot J'=(I_{3\times 3}-B)\cdot J,
\end{equation}
where $I_{3\times 3}$ is the $3$-dimensional identity operator.

This equation can be rewritten as
\begin{equation}\label{PF2}
D(t)\cdot J'(t)=Q(t)\cdot J(t),
\end{equation}
where
$$Q(t)=\begin{pmatrix}
  -\frac{1}{2}+32t & 9 & -10 \\
  0 & \frac{3}{2}+48t & -\frac{5}{2} \\
  0 & \frac{3}{2}-24t & -\frac{5}{2}+80t
  \end{pmatrix}\ \text{and} \ D(t)=64t^2-t.$$

Introducing new variables $X=t^{-\frac{1}{2}}J_1,\ Y=J_2,\ Z=J_3$ we have
\begin{equation}
\begin{cases}
\label{sys_XYZ}
D(t)\sqrt tX'&=9Y-10Z\\
D(t)Y'&=(\frac{3}{2}+48t)Y-\frac{5}{2}Z\\
D(t)Z'&=(\frac{3}{2}-24t)Y+(-\frac{5}{2}+80t)Z
\end{cases}
\end{equation}
\section{Proof of Theorem \ref{main_theorem}}
Take any Abelian integral $I(t)$. By Lemma \ref{lemma32} $I(t)=p_1(t)J_1(t)+p_2(t)J_2(t)+p_3(t)J_3(t)$ where $\deg p_i\leq \frac{n}{4}$.
Thus $I(t)=\sqrt tp_1X+p_2Y+p_3Z.$

Using \eqref{sys_XYZ}, we  obtain
\begin{equation}\label{frac_I}
\left(\frac{I}{p_1\sqrt t}\right)'=\frac{1}{Dp_1^2t\sqrt t}\cdot(\tilde p_1Y+\tilde p_2Z),
\end{equation}
where $\tilde p_i$ are some polynomials of degree less than or equal to $\frac{n}{2}+2$.

Recall that $Z=J_3=I_{0,3}(t)=\iint\limits_{\gamma(t)}y^3dx\wedge dy$, so $Z$ is positive for $t\in(0,\frac{1}{64})$ by Remark \ref{lem_pos}.
Hence the function $w=\frac{Y}{Z}$ is well-defined and by \eqref{sys_XYZ} satisfies the Riccati equation
\begin{equation}\label{ricati}
Dw'=\left(-\frac{3}{2}+24t\right)w^2+\left(4-32t\right)w-\frac{5}{2}.
\end{equation}

So for the function $S(t)=\tilde p_1w+\tilde p_2$ we have

$$DS'=\left(-\frac{3}{2}+24t\right)\tilde p_1w^2+(D\tilde p_1'+(4-32t)\tilde p_1)w+D\tilde p_2'-\frac{5}{2}\tilde p_1.$$

One can obtain
$$D\tilde p_1S'=\left(-\frac{3}{2}+24t\right)(S-\tilde p_2)^2+(D\tilde p_1'+(4-32t)\tilde p_1)(S-\tilde p_2)+\left(D\tilde p_2'-\frac{5}{2}\tilde p_1\right)\tilde p_1.$$

Thus the Riccati equation for the function $S(t)$ reads as
\begin{equation}\label{ricati_s}
D\tilde p_1S'=AS^2+BS+C,
\end{equation}
where $A,B$ and $C$ are polynomials and $\deg C\leq n+5.$
Now one can introduce new time $\tau$ and rewrite \eqref{ricati_s} as a system

\begin{equation}
\begin{cases}
\label{sys_ric}
\dot t&=D\tilde p_1\\
\dot S&=AS^2+BS+C,
\end{cases}
\end{equation}
where $\dot\varphi$ denotes $\frac{d\varphi}{d\tau}.$
Denote by $\Delta_j,$ $j=1,\dots,k\ (k\leq\deg\tilde p_1+1)$ the open intervals into which $(0,\frac{1}{64})$ is split by the zeros of $\tilde p_1$.
It is clear that in $\Delta_j$ between any two zeros of $S$ there is a zero of $C$. Let $\lambda_j$ be the number of zeros of $C$ on $\Delta_j$.
Thus the number of zeros of $S$ in $\Delta_j$ is less than or equal to $\lambda_j+1$.
So the number of zeros of $S$ in $(0,\frac{1}{64})$ is less than or equal to $\sum\limits_{j=1}^k(\lambda_j+1)\leq\deg C+\deg\tilde p_1+1.$
Thus it does not exceed $n+5+\frac{n}{2}+2+1=\frac{3}{2}n+8.$
By \eqref{frac_I} we obtain that on $(0,\frac{1}{64})$ the number of zeros of $\left(\frac{I}{p_1\sqrt t}\right)'$ does not exceed $\frac{3}{2}n+8.$

Denote by $\Xi_j,$ $j=1,\dots,l\ (l\leq\deg p_1+1)$ the open intervals into which $(0,\frac{1}{64})$ is split by the zeros of $p_1$.
It is clear that in $\Xi_j$ between any two zeros of $I$ (i.e. zeros of $\frac{I}{p_1\sqrt t}$) there is a zero of $\left(\frac{I}{p_1\sqrt t}\right)'$.
Let $l_j$ be the number of zeros of $\left(\frac{I}{p_1\sqrt t}\right)'$ on $\Xi_j$.
Thus the number of zeros of $I$ in $\Xi_j$ is less than or equal to $l_j+1$.
So the number of zeros of $I$ in $(0,\frac{1}{64})$ is less than or equal to $\sum\limits_{j=1}^l(l_j+1)$
and $\sum\limits_{j=1}^l l_j\leq \frac{3}{2}n+8.$
Thus the number of zeros of $I$ in $(0,\frac{1}{64})$ is less than or equal to $\frac{3}{2}n+8+l\leq \frac{3}{2}n+8+\frac{n}{4}+1=\frac{7}{4}n+9.$
This proves Theorem \ref{main_theorem}.

\end{document}